\documentclass[12pt,twoside]{article}
\usepackage{amsmath}
\usepackage{amsfonts}
\usepackage{amsthm}
\usepackage{amssymb}
\usepackage{graphicx}
\usepackage{multicol}

\begin{document}
\title{{\normalsize{\bf Ribbonness of Kervaire's sphere-link in homotopy 4-sphere and its consequences to 2-complexes}}}
\author{{\footnotesize Akio KAWAUCHI}\\
{\footnotesize{\it Osaka Central Advanced Mathematical Institute, Osaka Metropolitan University }}\\
{\footnotesize{\it Sugimoto, Sumiyoshi-ku, Osaka 558-8585, Japan}}\\
{\footnotesize{\it kawauchi@omu.ac.jp}}}
\date\, 
\maketitle
\vspace{0.25in}
\baselineskip=15pt
\newtheorem{Theorem}{Theorem}[section]
\newtheorem{Conjecture}[Theorem]{Conjecture}
\newtheorem{Lemma}[Theorem]{Lemma}
\newtheorem{Sublemma}[Theorem]{Sublemma}
\newtheorem{Proposition}[Theorem]{Proposition}
\newtheorem{Corollary}[Theorem]{Corollary}
\newtheorem{Claim}[Theorem]{Claim}
\newtheorem{Definition}[Theorem]{Definition}
\newtheorem{Example}[Theorem]{Example}

\begin{abstract} M. A. Kervaire showed that 
every group of deficiency $d$ and weight $d$ is the fundamental group of a smooth sphere-link of $d$ components in a smooth homotopy 4-sphere. In the use of the smooth unknotting conjecture and the smooth 4D Poincar{\'e} conjecture, any such sphere-link is shown to be  a sublink of a free ribbon sphere-link in the 4-sphere. Since every ribbon 
sphere-link in the 4-sphere is also shown to be a sublink of a free ribbon sphere-link in the 
4-sphere, Kervaire's sphere-link and the ribbon sphere-link are equivalent concepts. By applying this result to a ribbon disk-link in the 4-disk, it is shown that 
the compact complement of every ribbon disk-link in the 4-disk is aspherical. 
By this property, a ribbon disk-link presentation for every contractible 
finite 2-complex is introduced. By using this presentation, it  is shown that every connected subcomplex of a contractible finite 2-complex is aspherical (meaning partially yes for Whitehead aspherical conjecture).

\phantom{x}


\end{abstract}

\maketitle

\phantom{x}

\noindent{\bf 1. Introduction} 

 A group with finite presentation 
$<x_1, x_2, ..., x_n |\, r_1, r_2, ..., r_{n-d}>$ is called a group of  {\it deficiency} $d$. 
A group $G$ has {\it weight} $d$ if there are $d$ elements $w_1, w_2, ...,w_d$ in $G$
whose normal closure is equal to $G$, where the system of  elements $w_1, w_2, ...,w_d$ is  called a {\it weight system} of $G$. 
Let $X$ be a closed connected oriented smooth 4-manifold. A {\it sphere-link} or 
an $S^2$-{\it link} in $X$ is a disjoint sphere system smoothly embedded in $X$. 
A {\it surgery} of $X$ {\it along a loop system} $k_i\,(i=1,2,\dots, n)$ is the operation 
replacing a normal $D^3$-bundle system $k_i\times D^3\,(i=1,2,\dots, n)$ of 
$k_i\,(i=1,2,\dots, n)$ in $X$ by a normal $D^2$-bundle system 
$D^2_i\times S^2\,(i=1,2,\dots, n)$ of the 2-sphere system 
$K_i=0_i\times S^2\,(i=1,2,\dots, n)$ under the identifications that $\partial D^2_i =k_i\,(i=1,2,\dots, n)$ and $\partial D^3=S^2$. 
Let $X'$ be the smooth 4-manifold resulting from $X$ by this surgery. 
The spheres $K_i\,(i=1,2,\dots, n)$ form an $S^2$-link $K$ in $X'$. 
The 4-manifold $X'$ is said to be {\it obtained from the 4-manifold} $X$ 
{\it by surgery along a loop system} $k_i\,(i=1,2,\dots, n)$ in $X$, and conversely 
the 4-manifold $X$ is said 
to be {\it obtained from the 4-manifold} $X'$ {\it by surgery along a sphere system} $K$ 
in $X'$. Note that there are canonical fundamental group isomorphisms
\[\pi_1(X, v)\cong \pi_1(X\setminus k, v)\cong \pi_1(X'\setminus K, v)\] 
by general position. 
The {\it closed 4D handlebody of genus} $n$ is the 4-manifold 
\[Y^S = S^4 \#_{i=1}^n S^1 \times S^3_i\]
which is the connected sum of $S^4$ and $n$ copies 
$S^1 \times S^3_i\, (i=1,2,\dots,n)$ of the closed 4D handle $S^1\times S^3$. 
A {\it legged loop system with base point} $v$ in $X$ is a graph $\omega k$
of legged loops $\omega_i k_i\,(i=1,2,\dots,d)$ embedded in $X$ consisting of a disjoint simple loop system $k_i\,(i=1,2,\dots,d)$ and a leg system (=embedded path system) 
$\omega_i\,(i=1,2,\dots,d) $ such that 
$\omega_i$ connects the base point $v$ and a point $p_i \in k_i$ for every $i$ and the legs 
$\omega_i$ for all $i$ are made disjoint except for the base point $v$. 
The fundamental group $\pi_1(Y^S,v^S)$ is identified with the free group 
$<x_1,x_2,\dots,x_n>$ with basis $x_1,x_2,\dots,x_n$ represented by 
the {\it standard} legged loop system $\omega^S x$ 
of legged loops $\omega_i k_i\,(i=1,2,\dots,n)$ with base point $v^S$ in $Y^S$ 
using the standard loop $k_i=S^1\times \mathbf{1}_i$ of $S^1\times S^3_i$ 
and a leg $\omega_i$ joining $v^S$ and the point $(1, \mathbf{1}_i)\in 1 \times S^3_i$ 
not meeting $1 \times ( S^3_i\setminus\{ {\bf 1}_i \})$, for every $i$. 
A {smooth homotopy 4-sphere} is a smooth 4-manifold $M$ homotopy equivalent to the 4-sphere $S^4$. 
A {\it meridian system} of an $S^2$-link $K$ with $k$ components in $M$ is a legged loop system $\omega m$ with base point $v$ in $M\setminus L$ whose loop system $m$ consists of a meridian loop of every component of $K$. 
Kervaire showed the following theorem in \cite{Ker}\footnote{The condition that 
$H_1(G) = G/[G,G]$ is a free abelian group of rank $d$ is omitted since every  group $G$ of deficiency $d$ and weight $d$ has this condition.}.

\phantom{x} 

\noindent{\bf Kervaire's Theorem.} For every  group $G$ of 
deficiency $d$ and weight $d$, there is an $S^2$-link $K$ with $d$ components in a smooth homotopy 4-sphere $M$
such that there is an isomorphism $G\cong\pi_1(M\setminus K, v)$ sending the 
weight system to a meridian system of $K$.

\phantom{x}

The construction of an $S^2$-link in this theorem is explained as follows.

\phantom{x}

\noindent{\it Construction of Kervaire's $S^2$-link.} 
Let $<x_1, x_2, ..., x_n |\, r_1, r_2, ..., r_{n-d}>$ be a finite presentation of $G$ of deficiency $d$, and $w_1, w_2, ...,w_d$ a weight system of $G$. 
Let $G(n;n-d,d)$ be the {\it triple system} of the free group 
$<x_1, x_2, ..., x_n> $, the relator system $ r_1, r_2, ..., r_{n-d}$ written as words in $x_1, x_2, ..., x_n$ and a weight system $w_1, w_2, ...,w_d$ written as words in $x_1, x_2, ..., x_n$. 
Identify the free group $<x_1, x_2, ..., x_n>$ with the fundamental group 
$\pi_1(Y^S, v^S)$ of the 4D closed handlebody $Y^S$. 
Let $X$ be a 4-manifold obtained from $Y^S$ by surgery along a loop system 
$k(r_1), k(r_2), \dots, k(r_{n-d})$ in $Y^S$ representing the words 
$r_1, r_2, \dots, r_{n-d}$ in $\pi_1(Y^S, v^S)$. The fundamental group 
$\pi_1(X, v^S)$ has the presentation 
$<x_1, x_2, ..., x_n |\, r_1, r_2, ..., r_{n-d}>$ by Seifert-van Kampen theorem. 
Let $M$ be the 4-manifold obtained by surgery along a loop system 
$k(w_1),k(w_2), \dots, k(w_d)$ in $X$ representing 
the weight system $w_1, w_2, ...,w_d$ of $\pi_1(X, v^S)$. 
The manifold $M$ is a smooth homotopy 4-sphere by Seifert-van Kampen theorem. 
The $S^2$-link $K$ of $d$ components in $M$ is given by the core spheres   
$K_i =0_i\times\partial D^3\, (i=1,2,\dots, d)$  of  
$D^2_i\times \partial D^3$ replacing  $k(w_i)\times D^3\, (i=1,2,\dots, d)$.
The fundamental group $\pi_1(M\setminus K, v)$ is isomorphic to $\pi_1(X, v)\cong G$ 
by an isomorphism sending a meridian system of $K$ in $M$ to the weight system 
$w_1, w_2, ...,w_d$. 
This completes the construction of Kervaire's $S^2$-link. 

\phantom{x}

Kervaire's $S^2$-link $K$ in this construction is uniquely determined by the 
triple system $G(n;n-d,d)$ of the free group $<x_1, x_2, ..., x_n>$, 
the relator system $r_1, r_2, ..., r_{n-d}$ and the weight system $w_1, w_2, ...,w_d$, which 
is called Kervaire's $S^2$-link of {\it type} $G(n;n-d,d)$ or 
simply an $S^2$-link of {\it type} $G(n;n-d,d)$. 
For a smooth surface-link $L$ in $S^4$, the fundamental group $\pi_1(S^4 \setminus L, v)$ 
is a {\it meridian-based free group} if 
$\pi_1(S^4 \setminus L, v)$ is a free group with a basis represented by a meridian system 
of $L$ with base point $v$. A smooth surface-link $L$ in $S^4$ is a {\it trivial} surface-link if the components of $L$ bound disjoint handlebodies smoothly embedded in $S^4$. 
In this paper, Kervaire's $S^2$-link is studied by using  Smooth 4D Poincar{\'e} Conjecture and Smooth Unknotting Conjecture for a surface-link stated as follows:

\phantom{x}

\noindent{\bf Smooth 4D Poincar{\'e} Conjecture.}
Every 4D smooth homotopy 4-sphere $M$ is diffeomorphic to $S^4$. 

\phantom{x}

\noindent{\bf Smooth Unknotting Conjecture.} 
Every smooth surface-link $F$ in $S^4$ with a meridian-based free fundamental group 
$\pi_1(S^4 \setminus F, v)$ is a trivial surface-link.

\phantom{x}

The positive proofs of these conjectures are in \cite{K4} and \cite{K1,K2,K3}, respectively. 
{\it From now on, every smooth homotopy 4-sphere} $M$ {\it is identified with the 4-sphere} $S^4$. An $S^2$-link $L$ in $S^4$ is a {\it ribbon} $S^2$-link if $L$ is  an 
$S^2$-link obtained from a trivial $S^2$-link $O$ in $S^4$ by surgery along 
embedded 1-handles on $O$. See \cite[II]{KSS}, \cite{Yana1} 
for earlier concept of a ribbon surface-link. 
An $S^2$-link $L$ in $S^4$ is a {\it free} $S^2$-{\it link of rank} $n$ if 
the fundamental group $\pi_1( S^4 \setminus L, v)$ is a (not necessarily meridian based) free group of rank $n$. The following theorem is the first result of this paper.

\phantom{x}

\noindent{\bf Theorem~1.1.} The following three statements on an $S^2$-link $K$ with $d$ components in the 4-sphere $S^4$ are mutually equivalent:

\medskip

\noindent{(1)} The $S^2$-link $K$ is an $S^2$-link of type $G(n;n-d,d)$ for some $n$.

\medskip

\noindent{(2)} The $S^2$-link $K$ is a sublink with $d$ components of a free ribbon $S^2$-link 
of some rank $n$.

\medskip

\noindent{(3)} The $S^2$-link $K$ is a ribbon $S^2$-link with $d$ components. 

\phantom{x}

By combining Kervaire's Theorem with Theorems~1.1, the following characterization 
of the fundamental group $\pi_1(S^4 \setminus K, v)$ of a ribbon $S^2$-link $K$ in $S^4$ 
is obtained.

\phantom{x}

\noindent{\bf Corollary~1.2.} 
A group $G$ is a  group of deficiency $d$ and weight $d$
if and only if $G$ is isomorphic to 
the group $\pi_1(S^4\setminus K, v)$ of a ribbon $S^2$-link $K$ of $d$ components
in $S^4$ by an isomorphism  sending the weight system of $G$ to a meridian system of $K$. 

\phantom{x}

In the proof of Theorem~1.1, the claim that every free $S^2$-link is a free ribbon 
$S^2$-link is needed whose proof is done in \cite{K5}. For completeness of the present argument, 
this claim is moved to Appendix of this paper as {\it Free Ribbon Lemma} together with the proof. The proof of Theorem~1.1 is done in Section~2 by assuming Free Ribbon Lemma. 
A {\it trivial proper disk system} in the 4-disk $D^4$ is a disjoint proper disk system 
$D_i\, (i=1,2,\dots,n)$ in $D^4$ obtained by an interior push of a disjoint disk system $D^0_i\, (i=1,2,\dots,n)$ in the 3-sphere $S^3=\partial D^4$. 
A {\it ribbon disk-link} of $d$ components in $D^4$ is a disjoint proper disk system $L^D$ in $D^4$ 
of $d$ components which is obtained by an interior push of a disjoint disk system which is the union of a trivial proper disk system $D_i\, (i=1,2,\dots,n)$ in $D^4$ for some $n$ and a disjoint band system $b^0_j \, (j=1,2,\dots,n-d)$ in $S^3$ spanning the trivial link 
$\partial D_i\, (i=1,2,\dots,n)$ in $S^3$. The link $\partial L^D$ in $S^3$ is called a 
{\it classical ribbon link}.
By construction, the double of a ribbon disk-link $L^D$ of $k$ components in $D^4$ 
is a ribbon $S^2$-link $L$ of $k$ components in $S^4$. 
It is a standard fact that every ribbon $S^2$-link $(S^4,L)$ is considered as 
the double $(D^4\cup -D^4,L^D\cup -L^D)$ of a ribbon disk-link $(D^4,L^D)$ and its copy 
$(-D^4,-L^D)$, namely 
$(S^4,L)=(\partial(D^4\times I),\partial (L^D\times I)), I=[-1,1]$. 
To construct a ribbon disk-link $(D^4,L^D)$ from a ribbon $S^2$-link $(S^4,L)$, it is noted that 
a trivial $S^2$-link $O$ and embedded 1-handles to construct $L$ are always isotopically deformed   into a symmetric position with respect to the 
equatorial 3-sphere $S^3=\partial D^4=\partial(-D^4)$ in $S^4=D^4\cup -D^4$ 
 (see \cite[II]{KSS}). 
A {\it free ribbon disk-link of rank} $n$ is a ribbon disk-link $L^D$ in $D^4$ such that 
the fundamental group $\pi_1(D^4\setminus L^D,v)$ is a free group of rank $n$. 
In Lemma~3.1, it is shown that the inclusion $(D^4,L^D)\to (S^4,L)$ induces 
an isomorphism 
\[\pi_1(D^4\setminus L^D, v)\to \pi_1(S^4\setminus L, v).\] 
Thus, the $S^2$-link $L$ is a free ribbon $S^2$-link in $S^4$ if and only if 
the ribbon disk-link $L^D$ is a free ribbon disk-link in $D^4$. 
The {\it compact complement} of a ribbon disk-link $L^D$ in the 4-disk $D^4$ 
is the compact 4-manifold $E(L^D)=\mbox{cl}(D^4\setminus N(L^D))$ 
for a normal disk-bundle $N(L^D)=L^D\times D^2$ of $L^D$ in $D^4$. 
By Theorem~1.1, every ribbon $S^2$-link $K$ is 
a sublink of a free ribbon $S^2$-link $L$ of some rank $n$, so that every ribbon disk-link $K^D$ is a sublink of a free ribbon disk-link $L^D$ of some rank $n$ by Lemma~3.1. 
A connected complex is understood as a cell complex $P$ obtained from a 
bouquet of loops, called the {\it 1-skelton} $P^1$ of $P$, by adding $q(\geq 2)$-cells 
to $P^1$. A connected complex is {\it aspherical} if the universal covering space is contractible. 
A connected 2-complex $P$ is aspherical if and only if the second homotopy group 
$\pi_2(P,v)=0$. 
For a subcomplex $P'$ of a cell complex $P$, 
a {\it deformation retract} from $P$ to $P'$ 
is a map $r:P\to P'$ such that the composite map $ir:P\to P$ for the inclusion $i:P'\subset P$ 
is homotopic to the identity $1:P\to P$, where if the homotopy is further relative to $P'$, then 
the map $r$ is called a {\it strong deformation retract} from $P$ to $P'$ (see \cite{S}). 
It is shown in Lemma~3.2 that  for every free ribbon disk-link $L^D$ in $D^4$,  
there is a strong deformation retract 
\[r: E(L^D) \to \omega x\]
from the compact complement $E(L^D)$  to a legged $n$-loop system $\omega x$ in 
$E(L^D)$ representing the free group $\pi_1(E(L^D),v)=<x_1,x_2,\dots,x_n>$. 
Section~3 is devoted to explanations of Lemmas 3.1 and 3.2 on ribbon disk-links.
In Section~4, a decomposition of the 4-disk $D^4$ into a normal disk-bundle 
$N(L^D)=L^D\times D^2$ of a free ribbon disk-link $L^D$ and the compact complement $E(L^D)$ 
is considered. Let $Q(L^D)=E(L^D)\cup N(L^D)$ denote this decomposition of $D^4$. 
For a disk-link $L^D$ of $n$ components,  let $p_*=\{p_i|\,i=1,2,\dots,n\}$ be a set of 
$n$ points, one point taken from each component of $L^D$.
The strong deformation retract $N(L^D)\to p_*\times D^2$ shrinking $L^D$ into $p_*$ 
and the strong deformation 
retract $r: E(L^D) \to \omega x$ in Lemma~3.2 define a map 
\[\rho: Q(L^D) \to P(L^D)\]
with  $P(L^D)$ a finite 2-complex consisting of the 1-skelton $P(L^D)^1=\omega x$ 
and the 2-cells $p_*\times D^2$ attached by the attaching map 
$p_*\times \partial D^2\to\omega x$ defined by $r$. 
The map $\rho$ is called a {\it ribbon disk-link presentation} for the finite 2-complex $P(L^D)$. 
A {\it 1-full} subcomplex of a cell complex $P$ is a subcomplex $P'$ of $P$ such 
that the 1-skelton $(P' )^1$ of $P'$ is equal to the 1-skelton $P^1$ of $P$. 
For a sublink $K^D$ of $L^D$, let $N(K^D)=K^D\times D^2$ be the subbundle of the disk-bundle $N(L^D)$. 
Then the union $Q(K^D;L^D)=E(L^D)\cup N(K^D)$ is a decomposition of 
the compact complement $E(L^D\setminus K^D)$ 
of the sublink $L^D\setminus K^D$ of $L^D$ in $D^4$ and the ribbon disk-link presentation 
$\rho: Q(L^D) \to P(L^D)$ for $P(L^D)$ sends $Q(K^D;L^D)$ to a 1-full 2-subcomplex 
$P(K^D;L^D)$ of $P(L^D)$. 
Further, every 1-full 2-subcomplex of $P(L^D)$ is obtained from a sublink $K^D$ of $L^D$ in this way. The following theorem is shown in Section~4.

\phantom{x}

\noindent{\bf Theorem~1.3.} 
For every free ribbon disk-link $L^D$ in the 4-disk $D^4$, 
the ribbon disk-link presentation $\rho: Q(L^D) \to P(L^D)$ for the finite 2-complex 
$P(L^D)$  induces a homotopy equivalence 
$Q(K^D;L^D) \to P(K^D;L^D)$ for every sublink $K^D$ of  $L^D$  including 
$K^D=\emptyset$ and $K^D=L^D$. 
In particular, the finite 2-complex $P(L^D)$ is  contractible. 
Further, every contractible finite 2-complex $P$ is taken as $P=P(L^D)$ for 
a free ribbon disk-link $L^D$ in the 4-disk $D^4$. 

\phantom{x}

In Section~5, the following theorem is shown by using Theorem~1.3. 

\phantom{x}

\noindent{\bf Theorem~1.4.} The compact complement $E(K^D)$ of every 
ribbon disk-link $K^D$ in the 4-disk $D^4$ is aspherical.

\phantom{x}

The asphericity of the compact complement of a ribbon disk-knot in $D^4$ has been conjectured 
by Howie \cite{Howie} after having found some gaps on the arguments of Yanagawa \cite{Yana2} and Asano, Marumoto, Yanagawa \cite{AMY}. 
Since the fundamental group of an aspherical complex is torsion-free, the following corollary is obtained from Lemma~3.1 and Theorem~1.4. 

\phantom{x}

\noindent{\bf Corollary~1.5.} The fundamental group $\pi_1(S^4\setminus L, v)$ of every 
ribbon $S^2$-link in the 4-sphere $S^4$ is torsion-free.

\phantom{x}

This result gives the positive answer to the author's old question in \cite[II(pp.57-58)]{KSS}. 
The following corollary is obtained from Theorems~1.3 and 1.4, because if a 
connected subcomplex $P'$ of a contractible finite 2-complex $P$ is not 1-full, then 
a 1-full subcomplex $P''$ of $P$ is constructed from $P'$ by adding a bouquet of some loops in 
the 1-skelton $P^1$ of $P$ to $P'$, and $P''$ is aspherical if and only if $P'$ is aspherical. 

\phantom{x}

\noindent{\bf Corollary~1.6.} 
Every connected subcomplex of  every contractible finite 2-complex is aspherical.

\phantom{x}

This result is a partial positive confirmation of Whitehead aspherical conjecture \cite{White} 
claiming that every connected subcomplex of an aspherical 2-complex is aspherical. 

\phantom{x}

\noindent{\bf 2. Proof of Theorem~1.1}

The following lemma is a standard result obtained as a corollary of 
Smooth 4D Poincar{\'e} Conjecture and Smooth Unknotting Conjecture is shown in \cite[Corollary~1.5]{K4} without a mention of a legged loop system.

\phantom{x}

\noindent{\bf Lemma~2.1.} Every closed connected orientable smooth 4-manifold $Y$ with 
$\pi_1(Y, v)$ a free group and $H_2(Y; {\mathbf Z}) = 0$ is diffeomorphic to the closed 4D handlebody $Y^S$by a diffeomorphism $f:Y\to Y^S$ sending any given a legged loop system 
$\omega x$ 
with base point $v$ representing a basis $x_1,x_2,\dots,x_n$ of $\pi_1(Y, v)$  to  a 
standard legged loop system $\omega^S x$ of $Y^S$.
For any given spin structures on $Y$ and $Y^S$, the diffeomorphism $f$ can be taken 
spin-structure-preserving.

\phantom{x} 

\noindent{\bf Proof of Lemma~2.1.} 
(The proof   is moved  from \cite[Lemma~3.2]{K5} to here 
for completeness of the present argument.) 
Let $M$ be the 4-manifold obtained from $Y$ by surgery along the loop system $k(\omega x)$ 
of $\omega x$, which is identified with $S^4$ by Smooth 4D Poincar{\'e} Conjecture since it is a smooth homotopy 4-sphere by the van Kampen theorem and a homological argument. 
Let $L$ be the $S^2$-link in $S^4$ obtained from $k(\omega x)$ by the surgery. Then 
$\pi_1(S^4\setminus L,v)=<x_1,x_2,\dots,x_n>$ and the legged loop system $\omega x$ with 
base point $v$ in $Y$ is  a meridian system of $L$ in $S^4$ representing the basis 
$x_1,x_2,\dots,x_n$. 
By Smooth Unknotting Conjecture for an $S^2$-link, 
the $S^2$-link $L$ bounds disjoint 3-balls 
smoothly embedded in $S^4$ so that each 3-ball meets $\omega x$ with  just one 
transverse intersection  point in the loop system $k(\omega x)$ (see \cite{K3}). 
By the back surgery from $(M,L)$ to $(Y,k(\omega x))$, 
there is an orientation-preserving diffeomorphism 
$f: Y\to Y^S$ with $f(\omega x)=\omega^S x$. 
Given any spin structures on $Y$ and $Y^S$, note that there is 
an orientation-preserving spin-structure-changing diffeomorphism 
$:S^1\times S^3\to S^1\times S^3$ (see \cite{Gluck} for a
similar diffeomorphism on $S^1\times S^2$). 
Thus, by composing $f$ with the orientation-preserving spin-structure-changing diffeomorphisms 
on some connected summands of $Y^S$ which are copies of $S^1\times S^3$, the diffeomorphism 
$f: Y\to Y^S$ is modified into an orientation-preserving spin-structure-preserving diffeomorphism. 
This completes the proof of Lemma~2.1. $\square$

\phantom{x}

The proof of Theorem~1.1 is done as follows. 

\phantom{x}

\noindent{\bf 2.2: Proof of Theorem~1.1.}

\noindent{\it Proof of (1)$\to$(2).} Assume that an $S^2$-link $K$ of type $G(n;n-d,d)$ in $S^4$ for any $n$ is constructed from the triple system $G(n;n-d,d)$ consisting of the free basis 
$x_i\,(i=1,2,\dots, n)$, the relator system $r_i\,(i=1,2,\dots, n-d)$ written as words in 
$x_i\,(i=1,2,\dots, n)$ and a weight system $w_j\,(j=1,2,\dots, d)$ written as words in 
$x_i\,(i=1,2,\dots, n)$. 
The fundamental group $\pi_1(Y^S, v^S)$ of  $Y^S$ of rank $n$ is identified with 
the free group $<x_1, x_2, ..., x_n>$. 
Note that the elements $r_i, w_j\,(i=1,2,\dots,n-d; j=1,2,\dots,d)$ form a weight system of 
the free group $\pi_1(Y^S, v^S)$. 
Represent the elements $r_i, w_j\in \pi_1(Y^S, v^S)\,(i=1,2,\dots,n-d; j=1,2,\dots,d)$
by a disjoint simple loop system $k(r_i), k(w_j)\,(i=1,2,\dots,n-d; j=1,2,\dots,d)$ in $Y^S$. 
The 4-manifold $M$ obtained from $Y^S$ by surgery along the loop system 
$k(r_i), k(w_j)\,(i=1,2,\dots,n-d; j=1,2,\dots,d)$ is a smooth homotopy 4-sphere 
identified with $S^4$. 
Let $L $ be the $S^2$-link in $S^4$ of the sphere system 
$K(r_i), K(w_j)\,(i=1,2,\dots,n-d; j=1,2,\dots,d)$ occurring from the loop system 
$k(r_i), k(w_j)\,(i=1,2,\dots,n-d; j=1,2,\dots,d)$ by the surgery. 
The fundamental group $\pi_1(S^4\setminus L, v)$ is isomorphic to the free group 
$<x_1, x_2, ..., x_n>$ by an isomorphism sending a meridian system of $L$ to the
weight system $r_i, w_j\,(i=1,2,\dots,n-d; j=1,2,\dots,d)$. 
By Free Ribbon Lemma of Appendix, the $S^2$-link $L$ is a free ribbon $S^2$-link in $S^4$ 
of rank $n$. The sublink of $L$ consisting of the components $K(w_j)\,(j=1,2,\dots,d)$ is 
is just the $S^2$-link $K$ of type $G(n;n-d,d)$, which is  a sublink of 
the free ribbon $S^2$-link $L$ in $S^4$. 
This shows (1)$\to$(2). 

\noindent{\it Proof of (2)$\to$(1).} 
Let $K$ be a sublink of $d$ components of a free ribbon 
$S^2$-link $L$ of $n$ components in $S^4$ 
of rank $n$. Let $\pi_1(S^4\setminus L, v)= <x_1, x_2,\dots, x_n>$. 
Let $Y$ be the 4-manifold obtained from $S^4$ by surgery along $L$. By Lemma~2.1, $Y$ is identified with $Y^S$ of genus $n$ such that $\pi_1(S^4\setminus L, v)= <x_1, x_2,\dots, x_n>$ 
is identified with  $\pi_1(Y^S, v^S)$ by an isomorphism sending  a meridian system of 
$L$ in $S^4$ to a weight system of $\pi_1(Y^S, v^S)$. This means that 
the ribbon $S^2$-link $K$ is nothing but an $S^2$-link of type $G(n;n-d,d)$ for the 
triple system $G(n;n-d,d)$ consisting of the free group 
$\pi_1(Y^S, v)= <x_1, x_2,\dots, x_n>$, a relator 
system $r_1, r_2, ..., r_{n-d}$ coming from the meridian system of $L\setminus K$, and a 
weight system $w_1, w_2, ...,w_d$ coming from the meridian system of $K$. 
This shows (2)$\to$(1).

\noindent{\it Proof of (2)$\to$(3).} This proof is trivial. 

\noindent{\it Proof of (3)$\to$(2).} By definition, assume that a ribbon 
$S^2$-link $K$ of $d$ components in $S^4$ is obtained from a trivial $S^2$-link $O$ of 
$n$ components in $S^4$ by surgery along a 1-handle system $h$ on $O$. 
Let $O\times [0,1]$ be a collar of $O$ in $S^4$ where the 1-handle system $h$ meets  
only to  $O\times 0$, and 
$W = O\times [0,1] \cup h$ a $d$-component compact 3-manifold bounded by 
$K\cup O\times 1$. 
Let $K_i\,(i=1,2,\dots,d)$ be the components of $K$. Let $O'$ be a 
sublink of $O\times 1$ of $n-d$ components obtained by removing any one component 
of $O\times 1$ from the boundary of the component of $W$ containing the component $K_i$ for every $i$. Then there are isomorphisms 
\[\pi_1((S^4 \setminus W, v) \to \pi_1(S^4 \setminus K\cup O'), v)
\quad\mbox{and}\quad \pi_1(S^4 \setminus W, v)\to \pi_1(S^4\setminus O, v). \]
This is because there are deformation retracts from $W$ to a 2-complex 
consisting of $K\cup O'$ and some spanning arcs and from $W$ to a 2-complex consisting of $O$ and some spanning arcs, and the spanning arcs do not affect the fundamental group. 
Since $\pi_1(S^4\setminus O, v).$ is a free group of rank $n$,  the $S^2$-link $L=K\cup O'$ 
of $n$ components
is a free ribbon $S^2$-link of rank $n$ in $S^4$ containing $K$ as a sublink. 
This shows (3)$\to$(2). 

This completes the proof of Theorem~1.1. $\square$

\phantom{x}

\noindent{\bf 3. Basic Lemmas of ribbon disk-links}

For a ribbon disk-link $(D^4, L^D)$ of a ribbon $S^2$-link $(S^4,L)$, 
let $\alpha$ be the reflection of $(S^4,L)$ exchanging $(D^4, L^D)$ and the other copy 
$(-D^4,-L^D)$ in $(S^4,L)$. 
Although the following lemma may be more or less known (cf. \cite{Yana2}), the proof is given here for convenience. 

\phantom{x}

\noindent{\bf Lemma~3.1.} For a ribbon disk-link $L^D$ in $D^4$ of a ribbon $S^2$-link $L$ in $S^4$, the inclusion $(D^4, L^D)\to (S^4,L)$ induces an isomorphism
\[\pi_1(D^4\setminus L^D, v) \to \pi_1( S^4\setminus L, v).\]

\phantom{x}

\noindent{\bf Proof of Lemma~3.1.} Use the retraction $S^4\setminus L \to D^4\setminus L^D$ induced from the quotient by the reflection $\alpha$. Then the canonical homomorphism 
$\pi_1(D^4\setminus L^D, v) \to \pi_1( S^4\setminus L, v)$ is shown to be a monomorphism. 
On the other hand, for the copy $(-D^4,-L^D)$ of $(D^4, L^D)$, 
the inclusion $(\partial (-D^4), \partial (-L^D))\to (-D^4, -L^D)$ induces an epimorphism 
$\pi_1(\partial (-D^4)\setminus \partial (-L^D),v)\to \pi_1(-D^4\setminus -L^D,v)$ 
by the definition of ribbon disk-link and Seifert-van Kampen theorem.
This means that the canonical monomorphism 
$\pi_1(D^4\setminus L^D, v) \to \pi_1( S^4\setminus L, v)$ is also an epimorphism and thus, 
an isomorphism. $\square$

\phantom{x}

The {\it 4D handlebody of genus} $n$ is the 4-manifold 
\[Y^D = D^4 {}_{\partial}\#_{i=1}^n S^1 \times D^3_i\]
which is the boundary connected sum of $D^4$ and $n$ copies 
$S^1 \times D^3_i\, (i=1,2,\dots,n)$ of the 4D handle $S^1\times D^3$. 
By using the asphericity of $Y^D$, the following lemma is obtained.

\phantom{x}

\noindent{\bf Lemma~3.2.}
For every free ribbon disk-link $L^D$ of rank $n$ in $D^4$, 
there is a strong deformation retract
\[r: E(L^D)\to \omega x\]
from the compact complement $E(L^D)$ to a legged $n$-loop system $\omega x$ 
with base point $v$ in  $E(L^D)$ 
representing any basis $ x_1,x_2,\dots,x_n$ of the free group $\pi_1(E(L^D),v)$.

\phantom{x}

\noindent{\bf Proof of Lemma~3.2.} Let $L$ be the free ribbon $S^2$-link of rank $n$ in $S^4$ 
obtained by taking the double of $(D^4,L^D)$.
Note that the double 
\[Y=\partial (E(L^D)\times I)=E(L^D)\times \{-1\}\cup (\partial E(L^D))\times I 
\cup E(L^D)\times \{1\}\]
of $E(L^D)$ is diffeomorphic to the 4-manifold $Y'$ obtained from $S^4$ by surgery along $L$. 
Since there is a canonical isomorphism 
$\pi_1(S^4\setminus L,v)=<x_1,x_2,\dots,x_n>\to \pi_1(Y',v)$ 
and $H_2(Y';{\mathbf Z})=0$, the 4-manifold $Y'$ is identified with $Y^S$ under the canonical identities $\pi_1(E(L^D,v)= \pi_1(Y^S,v)=<x_1,x_2,\dots,x_n>$ by Lemmas~2.1 and 3.1. 
Let $\omega x$ be a legged $n$-loop system in $E(L^D)$, and $-\omega x$ a copy of 
$\omega x$ in the copy $-E(L^D)$ of $E(L^D)$ in $Y'=Y^S$. Note that $\pm \omega x$ are isotopically deformed into the standard $n$-loop system in $Y^S$. 
Let $N(\omega x)$ be a regular neighborhood 
of $\omega x$ in $E(L^D)$, and $N(-\omega x)$ the copy of $N(\omega x)$ 
in the copy $-E(L^D)$. Since $N(\omega x)$ 
is diffeomorphic to the 4D handlebody $Y^D$ of genus $n$, it is shown that 
the compact complement $E(L^D)^+=\mbox{cl}(Y^S\setminus N(-\omega x))$ 
is  diffeomorphic to $Y^D$ and 
the compact complement $H=\mbox{cl}(Y^S\setminus N(\omega x)\cup N(-\omega x))$ 
is diffeomorphic to the product $Z^S\times I$ for the closed 3D handlebody 
$Z^S = S^3 \#_{i=1}^n S^1 \times S^2_i$ of genus $n$. 
Note that the reflection $\alpha$ in $Y^S$ exchanging $E(L^D)$ 
and $-E(L^D)$ induces a reflection in $H$ whose fixed point set is the boundary 
$Z(\partial L^D)=\partial E(L^D)$ of $E(L^D)$.  
Let $H'$ be one of the two  3-manifolds  obtained  by splitting $H$ along $Z(\partial L^D)$ 
such that $E(L^D)^+=E(L^D)\cup H'$.  Then  $H=H'\cup\alpha(H')$. 
By \cite{K0}, the 3-manifold $Z(\partial L^D)$ is an imitation of $Z^S$ which has 
the property that the inclusion homomorphism $\pi_1(Z^S,v)\to \pi_1(H',v)$ is an isomorphism 
and any covering triad $(\tilde H'; \tilde Z(\partial L^D), \tilde Z^S)$ of the triad 
$(\tilde H'; Z(\partial L^D),Z^S)$ is a homology cobordism. This means that 
the inclusion $i:E(L^D)\to E(L^D)^+$ is a homotopy equivalence by Seifert-van Kampen theorem and the universal covering lift 
$\tilde i: \tilde E(L^D)\to \tilde E(L^D)^+$ induces an isomorphism 
$\tilde i_*: H_*(\tilde E(L^D);{\mathbf Z})\to H_*(\tilde E(L^D)^+;{\mathbf Z})$ because 
\[H_*(\tilde E(L^D)^+, \tilde E(L^D);{\mathbf Z})\cong 
H_*(\tilde H',Z(\partial L^D);{\mathbf Z})=0\]
by the excision isomorphism. Thus, $E(L^D)$ is homotopy equivalent to the 
legged $n$-loop system $\omega x$.
For a polyhedral pair $(P,P')$, if the inclusion $i:P'\subset P$ is homotopy equivalent, then 
there is a strong deformation retract $r:P\to P'$ (see \cite[p.~31]{S}). 
Thus, there is a strong deformation retract $r:E(L^D)\to \omega x$.
$\square$

\phantom{x}

In Lemma~3.2, note that in general the compact complement $E(L^D)$ of a free ribbon disk-link $L^D$ in $D^4$ is not diffeomorphic to $Y^D$. For example, the Kinoshita-Terasaka knot $k_{KT}$ in $S^3$ bounds a free ribbon-disk knot $K^D$ of rank one in $D^4$. 
Since the 3-manifold $Z(\partial K^D)$ which is the $0$-surgery manifold 
of $k_{KT}$  is not diffeomorphic to $Z^S=S^1\times S^2$ by the solution of 
property R conjecture (see \cite{Gabai}),  
the compact complement $E(K^D)$ is not diffeomorphic to $Y^D$  (see \cite{K0}).  

\phantom{x}

\noindent{\bf 4. Proof of Theorem~1.3}

The proof of Theorem~1.3 is done as follows. 

\phantom{x}

\noindent{\bf 4.1: Proof of Theorem~1.3.} 
Identifications 
\[\pi_1(E(L^D),v)=\pi_1(\omega x)=<x_1,x_2,\dots,x_n>\] 
are fixed by the strong deformation retract $r:E(L^D),v)\to \omega x$. 
The ribbon disk-link presentation $\rho : Q(L^D) \to P(L^D)$ for $P(L^D)$ induces an isomorphism 
$\rho_{\#}:\pi_1(Q(K^D;L^D),v) \to \pi_1(P(K^D;L^D),v)$ for every sublink $K^D$ of $L^D$ 
including $K^D=\emptyset$ and $K^D=L^D$ by Seifert-van Kampen theorem, because 
the strong deformation retract $r:E(L^D),v)\to \omega x$ induces 
the identical word system $r_*=\{r_1,r_2,\dots,r_n\}$ of the loop system 
$p_*\times \partial D^2$ in 
$<x_1,x_2,\dots,x_n>$ by the attaching map $r: p_*\times \partial D^2\to\omega x$ of 
the 2-cell system $p_*\times D^2$.
 In particular, $\pi_1(P(L^D),v)=<x_1,x_2,\dots,x_n|\,r_1,r_2,\dots,r_n>=\{1\}$. 
Let $\tilde \rho : \tilde Q(K^D;L^D) \to \tilde P(K^D;L^D)$ be the universal covering lift 
of $\rho : Q(K^D;L^D) \to P(K^D;L^D)$. By Mayer-Vietoris homology sequence, 
$H_m(\tilde Q(K^D;L^D);{\mathbf Z})=0$ for all $m\geq 3$ 
and $\tilde \rho$ induces an isomorphism 
$\tilde \rho_*: H_2(\tilde Q(K^D;L^D);{\mathbf Z})\to H_2(\tilde P(K^D;L^D);{\mathbf Z})$ 
for every sublink $K^D$ of $L^D$ 
including $K^D=\emptyset$ and $K^D=L^D$.
Thus, $\rho : Q(K^D;L^D) \to P(K^D;L^D)$ is a homotopy equivalence for every sublink $K^D$ of $L^D$ 
including $K^D=\emptyset$ and $K^D=L^D$. In particular, $P(L^D)$ is a finite contractible 2-complex. 
Let $P$ be a contractible finite 2-complex obtained from the 1-skelton $P^1=\omega x$, a legged $n$ loop system  with base point $v$, so that  $\pi_1(P^1, v)=<x_1, x_2,\dots., x_n>$. 
Assume that $P$ is obtained from $P^1$ by attaching 2-cells $e_1, e_2,\dots, e_n$. 
Since $\pi_1(P, v) ={1}$, the 2-complex $P$ provides the
triple system $G(n;0,n)$ in the construction of Kervaire's 2-link which consists of the free group 
$<x_1, x_2,\dots., x_n>$, the empty relator set and the weight system $w_1, w_2, ...,w_n$ 
given by the attaching data of $e_1,e_2,\dots, e_n$ to $P^1$. 
By Theorem~1.1, 
there is a free ribbon $S^2$-link $(S^4,L)$ with an isomorphism 
$\pi_1(S^4\setminus L, v) \cong <x_1, x_2,\dots, x_n>$ 
sending a meridian system of $L$ to the weight system $w_1, w_2, ...,w_n$. 
By Lemma~3.1, there is a free ribbon disk-link $(D^4,L^D)$ with an isomorphism 
$\pi_1(D^4\setminus L^D, v) \cong <x_1, x_2,\dots, x_n>$ 
sending a meridian system of $L^D$ in $D^4$ to the weight system $w_1, w_2, ...,w_n$. 
By Lemma~3.2, there is a strong deformation retract $r:E(L^D)\to P^1=\omega x$, which 
induces a ribbon-disk presentation $\rho:Q(L^D)\to P(L^D)$ for $P(L^D)=P$ because the loop 
system $p_*\times \partial D^2$ is just the meridian system of $L^D$. 
$\square$

\phantom{x}

\noindent{\bf 5. Proof of Theorem~1.4}

The proof of Theorem~1.4 is done as follows.

\phantom{x}

\noindent{\bf 5.1: Proof of Theorem~1.4.} 
Let $K^D$ be a ribbon disk-link in $D^4$ of $d$ components, and 
$S(*)$ any immersed 2-sphere in $E(K^D)$. 
It suffices to show that there is a free ribbon disk-link $L^D$ in $D^4$ of some rank $n$ which contains $K^D$ as a sublink and is disjoint from $S(*)$. 
This is because $S(*)\subset E(L^D)\subset E(K^D)$ meaning that 
$S(*)$ is null-homotopic in $E(L^D)$ and hence in $E(K^D)$ 
since  $\pi_2(E(L^D),v)=0$ by Lemma~3.2, so that $\pi_2(E(K^D),v)=0$ 
meaning that $E(K^D)$ is aspherical, for  
$E(K^D)$ is homotopy equivalent to a 2-complex by Theorem~1.3.

The pair $(D^4,S^3)$ is considered as the one-point compactification 
of the pair 
$({\mathbf R}^3[0,+\infty),{\mathbf R}^3)$ of the upper-half 4-space
\[{\mathbf R}^3[0,+\infty)=\{(x_1,x_2,x_3,t)|\, -\infty<x_i<+\infty\, (i=1,2,3), t\geq 0\}\] 
and the 3-space 
\[{\mathbf R}^3=\{(x_1,x_2,x_3)|\, -\infty<x_i<+\infty\, (i=1,2,3)\}. \] 
Also, $K^D$ and $S(*)$ are considered in ${\mathbf R}^3[0,+\infty)$. 
By the motion picture method \cite[I]{KSS}, assume that a normal form of 
the disk-link $K^D$ in 
$({\mathbf R}^3[0,+\infty)$ is given as follows: 
\[
K^D\cap {\mathbf R}^3[t]=\left\{
\begin{array}{rl}
\emptyset,&\quad \mbox{for $t>2$},\\
d_*[t],& \quad \mbox{for $t=2$},\\
o_*[t],& \quad \mbox{for $1<t<2$},\\
(o_*\cup b_*)[t],& \quad \mbox{for $t=1$},\\
k^D[t],& \quad \mbox{for $0\leq t<1$},
\end{array}\right.
\]
where $d_*$ is a disjoint trivial disk system of $m$ disks $d_i\, (i=1,2,\dots,m)$ for some $m$ 
in ${\mathbf R}^3$ with $o_*=\partial d_*$, 
$b_*$ is a disjoint band system of $m-d$ bands $b_j \, (j=1,2,\dots,m-d)$ 
in ${\mathbf R}^3$ spanning the trivial loop system $o_*$ used for a fusion operation, 
and $k^D$ is a ribbon link in ${\mathbf R}^3$ of $d$-components obtained from $o_*$ 
by surgery along the band system $b_*$ as a fusion. 
By the proof of Theorem~1.1 and Lemma~3.1, there is a free ribbon disk-link $L^D$ in 
${\mathbf R}^3[0,+\infty)$ of some rank $n$ 
such that $L^D=K^D\cup C^D$ for a trivial disk system $C^D$ in 
${\mathbf R}^3[0,+\infty)$ 
whose normal form  is given as follows by extending the normal form of $K^D$:
\[
L^D\cap {\mathbf R}^3[t]=\left\{
\begin{array}{rl}
\emptyset,&\quad \mbox{for $t>2$},\\
(d_*\cup d^C)[t],& \quad \mbox{for $t=2$},\\
(o_*\cup o^C)[t],& \quad \mbox{for $1<t<2$},\\
(o_*\cup b_*\cup o^C)[t],& \quad \mbox{for $t=1$},\\
(k^D\cup o^C)[t],& \quad \mbox{for $0\leq t<1$},
\end{array}\right.
\]
where $d^C$ is a disjoint disk system in ${\mathbf R}^3$ with 
$o^C=\partial d^C$. Note that the disk systems $d_*$ and $d^C$ are disjoint, but 
in general the band system $b_*$  meets the interior of $d^C$ in a disjoint arc system. 
By pulling down a neighborhood of every double point of $S(*)$ into 
${\mathbf R}^3[0]$, the immersed 2-sphere $S(*)$ is changed into a 
non-immersed singular 2-sphere in ${\mathbf R}^3[0,+\infty)$, but 
a normal form of the union $K^D\cup  S(*)$ in ${\mathbf R}^3[0,+\infty)$ extending the normal form of $K^D$  is given  as follows (see  \cite[I]{KSS}): 
\[
(K^D\cup S(*))\cap {\mathbf R}^3[t]=\left\{
\begin{array}{rl}
\emptyset,&\quad \mbox{for $t>2$},\\
(d_*\cup d^{S(*)})[t],& \quad \mbox{for $t=2$},\\
(o_*\cup o^{S(*)})[t],& \quad \mbox{for $1<t<2$},\\
(o_*\cup b_*\cup c^{S(*)}\cup b^{S(*)})[t],& \quad \mbox{for $t=1$},\\
(k^D\cup c^{S(*)})[t],& \quad \mbox{for $0<t<1$},\\
(k^D\cup e^{S(*)})[t],& \quad \mbox{for $t=0$},
\end{array}\right.
\]
where $d^{S(*)}$ is a disjoint band system in ${\mathbf R}^3$ with 
$o^{S(*)}=\partial d^{S(*)}$, $b^{S(*)}$ is a disjoint band system spanning 
$o^{S(*)}$ in ${\mathbf R}^3$, $c^{S(*)}$ is  a split union of a split Hopf link system $c^{H(*)}$ and a trivial link system $c^{o(*)}$ in ${\mathbf R}^3$ obtained from $o^{S(*)}$ by surgery along $b^{S(*)}$, and 
$e^{S(*)}$ is a split union of a disjoint Hopf disk pair system
bounded by $c^{H(*)}$ and a disjoint disk system bounded by $c^{o(*)}$ in 
${\mathbf R}^3$, where a {\it Hopf disk pair} means a disk pair with a clasp singularity  in ${\mathbf R}^3$ bounded by a Hopf link. By construction, note that $e^{S(*)}$ is split from $k^D$.  
By an isotopic move of the union of the disk system $d^C$ and a neighborhood of 
the arc system 
$b_*\cap d^C$ in $b_*$ in ${\mathbf R}^3$ keeping the disk system $d_*$ fixed, 
it can be assumed that 
\[d^C\cap(d_*\cup e^{S(*)}\cup b^{S(*)})=\emptyset.\] 
Then the link $o_*\cup o^{S(*)}\cup o^C$ is a trivial link in ${\mathbf R}^3$. 
In general the disk system $d^C$ meets the interior of the disk system $d^{S(*)}$.  However, by Horibe-Yanagawa lemma in \cite[I]{KSS}, 
even if the disk systems $d_*, d^{S(*)}, d^C$ are replaced by any disjoint disk systems bounded by the trivial link $o_*\cup o^{S(*)}\cup o^C$ in 
${\mathbf R}^3$, the union $K^D\cup S(*)$ and the free ribbon disk-link $L^D$ in 
 do not change up to ambient isotopies (with compact supports) of 
${\mathbf R}^3[0,+\infty)$  keeping ${\mathbf R}^3[0]$  fixed. 
This means that the disjoint union $K^D\cup S(*)$ extends to a disjoint union $L^D\cup S(*)$ for a free ribbon disk-link $L^D$, so that 
$S(*)\subset E(L^D)\subset E(K^D)$, and thus, $E(K^D)$ is aspherical. This completes the proof of Theorem~1.4. 
$\square$

\phantom{x}

\noindent{\bf Appendix: Free Ribbon Lemma}

The purpose of this appendix is to prove the following lemma.

\phantom{x}

\noindent{\bf Free Ribbon Lemma.} 
Every free $S^2$-link $L$ in $S^4$ is a ribbon $S^2$-link.

\phantom{x} 

\noindent{\bf Proof of Free Ribbon Lemma.} 
The following observation is used to determine a ribbon $S^2$-link. 

\phantom{x}

\noindent{\bf (A.1)} Let $(S^3_i)^{(1+m_i)}\, (i=1,2,\dots,n)$ be a 
system of mutually disjoint compact $(1+m_i)$-punctured 3-spheres in $S^4$ 
such that the boundary $\partial (S^3_i)^{(1+m_i)}$ is 
the union of the  component $K_i$ and an $S^2$-link $O_i$ of $m_i$ components. If the union $O=\cup_{i=1}^n O_i$ is a trivial $S^2$-link in $S^4$, then 
the $S^2$-link $L=\cup_{i=1}^n K_i$ is a ribbon $S^2$-link in $S^4$. 

\phantom{x}

\noindent{\bf Proof of (A.1).} Let $K'_i$ be a 2-sphere obtained from $O_i$ by 
surgery along mutually disjoint 1-handles $h_i\, (i=1,2,\dots, m_i-1)$ in $(S^3_i)^{(1+m_i)}$, 
whose closed complement is diffeomorphic to the spherical shell $S^2\times[0,1]$.
This means that the component $K_i$ with reversed orientation is isotopic to the 
2-sphere $K'_i$ in $(S^3_i)^{(1+m_i)}$. 
This shows that $L=\cup_{i=1}^n K_i$ is a ribbon $S^2$-link in $S^4$, completing the proof 
of (A.1). $\square$

\phantom{x}

Let $K_i\,(i=1,2,\dots,n)$ be the components of a free $S^2$-link $L$ in $S^4$. 
Let $Y$ be the 4-manifold obtained from $S^4$ by surgery along $L$. 
Let $k_i\,(i=1,2,\dots,n)$ be the loop system in $Y$ produced from  
$K_i\,(i=1,2,\dots,n)$ by the surgery. Since the fundamental group $\pi_1(Y, v)$ 
is a free group and $H_2(Y; {\mathbf Z}) = 0$, the 4-manifold $Y$ is identified with 
$Y^S$ by Lemma~2.1.  The 3-sphere $1\times S^3_i$ of the connected 
summand $S^1\times S^3_i$ of $Y^S$ is fixed and denoted by $S^3_i$. 
Let $x_i\, (i=1,2,\dots, n)$ be the basis of $\pi_1(Y^S,v) $ represented by a standard legged loop system $\omega^S x$ with vertex $v=v^S$. 
Let $k(\omega^S x)=\{ k^S_i|\, i=1,2,\dots, n\}$ be the loop system of $\omega^S x$.
Let $\omega m=\{\omega_i m_i|\,i=1,2,\dots, n\}$ be a  meridian system  
with vertex $v$ of the components $K_i\,(i=1,2,\dots,n)$ of $L$ in $S^4$. 
The meridian system $\omega m$ is taken in $Y^S$ as a legged loop system with  
loop system  $k(\omega m)=\{m_j|\, j=1,2,\dots,n\}$ parallel  to the loop system 
$k_i\,(i=1,2,\dots,n)$ in $Y^S$. 
Assume that the meridian system $\omega m$ in $Y^S$  is made disjoint from 
$\omega x$ except for the vertex $v$ and meets $S^3_i\,(i=1,2,\dots,n)$ only in the loop system $k(\omega m)$ transversely. 
Let $y_i\,(i=1,2,\dots,n)$ be the elements of $\pi_1(Y^S,v)$ represented by 
$\omega_i m_i\, (i=1,2,\dots, n)$. 
By Nielsen transformations of the basis $x_i\,(i=1,2,\dots,n)$, assume that 
the product $x^{-1}_i y_i$ is in the commutator subgroup $[\pi_1(Y^S,v),\pi_1(Y^S,v)]$ of 
$\pi_1(Y^S,v)$ for every $i$ (see \cite{MKS}). 
For the 3-sphere $S^3_i$, consider all the loops $m_j$ with $m_j\cap S^3_i\ne\emptyset$. 
For a point $p\in m_j\cap S^3_t\,(t\ne i)$, let $I(p)$ be an arc neighborhood of $p$ in a parallel $k^S_t(p)$ of $k^S_t$ and then replace the arc $I(p)$ with the arc 
$\mbox{cl}(k^S_t(p)\setminus I(p))$. 
Let $\tilde m_j$ be a loop obtained from $m_j$  by doing this operation on $m_j$ for every 
$t\, (t\ne i)$ and every point $p\in m_j\cap S^3_t$. 
For every $i\,(i=1,2,\dots,n)$, let $m(S^3_i)$ be the system of the loops $\tilde m_j$ in $Y^S$ obtained from all the loops $m_j$ with $m_j\cap S^3_i\ne \emptyset$, where the loops $m_j$ with $m_j\cap S^3_i= \emptyset$ are discarded. 
There is a smoothly embedded annulus $A_i$ with $\partial A_i=(-k^S_i)\cup \tilde m_i$ 
in the open 4-manifold 
\[ Y^S_{)i(}=Y^S\setminus \cup_{1\leq t(\ne i)\leq n} S^3_t\]
because the fundamental group $\pi_1(Y^S_{)i(},v)$ is an infinite cyclic group and the loop 
$\tilde m_i$ is homotopic to $k^S_i$ in $Y^S_{)i(}$. 
The annulus $A_i$ meets $S^3_i$ transversely with disjoint simple loops and simple arcs. 
Let $\alpha_{is}\, (s= 1,2,\dots, n_i)$ be the arc system of the intersection 
$A_i\cap S^3_i$ where $\alpha_{i1}$ joins the point $p^S_i=k^S_i\cap S^3_i$ to a point of 
the loop $\tilde m_i$ and the arc $\alpha_{is}$ with $s>1$ joins two points of $\tilde m_i$. 
For $j$ with $j\ne i$, the loop $\tilde m_j$ is null-homotopic in $Y^S_{)i(}$ 
and hence bounds a disk 
$D_{j i}$ in $ Y_{)i(}$ which meets $S^3_i$ transversely with disjoint simple loops and simple arcs. Let $\alpha_{j i s}\, (s=1,2,\dots, n_{j i})$ be the arc system of the intersection 
$D_{j i}\cap S^3_i$ each of which joints two points of $\tilde m_j$. 
The annulus $A_i$ and the disk $D_{j i}$ with $i\ne j$ are made disjoint while fixing 
the intersection with $S^3_i$ in $Y^S$ for all $i,j$ 
by  doing double point cancellations using free boundary arcs 
while fixing the intersection with $S^3_i$ for $m_j$. 
The following observation helps clarify the relationship 
between the point system $m(S^3_i)\cap S^3_i$ 
and the arc system $(A_i\cup D_{j i})\cap S^3_i$ for all $j$ with $j\ne i$. 

\phantom{x}

\noindent{\bf Observation (A.2)} 
Let  $\partial \alpha_{is}=\{q_s,q'_s\}\, (s= 1,2,\dots,n_i)$ with $q_1=p^S_i$ for the arc system 
$\alpha_{is}\, (s= 1,2,\dots, n_i)$ of $A_i\cap S^3_i$.
Then the  open arc of $\tilde m_i$ that is separated by any couple $\{q_s,q'_s\}$ with $s>1$ and does not contain the point $q'_1$ meets $S^3_i$ with intersection number $0$. 
Conversely, let $\{q_s,q'_s\}\, (s=1,2,\dots,n_i)$ be any system of couples of distinct points 
with $q_1=p^S_i$ such that the union of these points matches 
the set $(k^S_i\cup \tilde m_i)\cap S^3_i$  and  
the open arc of $\tilde m_i$ that is divided by any couple $\{q_s,q'_s\}$ with $s>1$ and does not contain the point $q'_1$ meets $S^3_i$ with intersection number $0$. 
Then $\{q_s,q'_s\}\, (s=1,2,\dots,n_i)$ is realized by 
$\partial \alpha_{is}=\{q_s,q'_s\}\, (s= 1,2,\dots,n_i)$ of the arc system 
$\alpha_{is}\, (s=1,2,\dots, n_i)$ 
of $A_i\cap S^3_i$ for an annulus $A_i$ with $\partial A_i=(-k^S_i)\cup \tilde m_i$ in 
$Y^S_{)i(}$. 
Let $\partial \alpha_{j is}=\{q_s,q'_s\}\, (s= 1,2,\dots,n_{j i})$ for the arc system 
$\alpha_{j is}\, (s=1,2,\dots, n_{j i})$ of $D_{j i}\cap S^3_i$. 
Then every open arc of $\tilde m_j$ divided by any couple $\{q_s,q'_s\}$ 
meets $S^3_i$ with intersection number $0$. 
Conversely, let $\{q_s,q'_s\}\, (s=1,2,\dots,n_{j i})$ be any system of couples of distinct points  
such that the union of these points matches the set $\tilde m_j\cap S^3_i$ and   
every open arc of $\tilde m_j$ which is divided by any couple $\{q_s,q'_s\}$ meets $S^3_i$ with intersection number $0$. 
Then $\{q_s,q'_s\}\, (s=1,2,\dots,n_{j i})$ is realized by 
$\partial \alpha_{j is}=\{q_s,q'_s\}\, (s= 1,2,\dots,n_{j i})$
of the arc system $\alpha_{j is}\, (s=1,2,\dots, n_{j i})$ 
of $D_{j i}\cap S^3_i$ for a disk $D_{j i}$  with $\partial D_{j i}=\tilde m_j$ in $ Y^S_{)i(}$.

\phantom{x}

Let $B(\alpha_{is})\, (s= 1,2,\dots, n_i)$ be disjoint 3-ball neighborhoods of the arcs 
$\alpha_{is}\, (s= 1,2,\dots, n_i)$ in $S^3_i$, and 
$B(\alpha_{j is})\, (s=1,2,\dots, n_{j i})$ disjoint 3-ball neighborhoods of the arcs
$\alpha_{j is}\, (s=1,2,\dots, n_{j i})$ in $S^3_i$. 
Let $S(\alpha_{is})=\partial B(\alpha_{is})\,(s= 1,2,\dots, n_i)$ and 
$S(\alpha_{j is})=\partial B(\alpha_{j is})\, (s=1,2,\dots, n_{j i})$ be 
the boundary 2-spheres of them. 
The $S^2$-link $L$ in $S^4$ with meridian system $\omega m$  is recovered from $Y^S$ 
by the back surgery along the loop system $k_i\,(i=1,2,\dots,n)$ in $Y^S$.  
Since the 2-spheres $S(\alpha_{is})\,(s= 1,2,\dots, n_i)$ and 
$S(\alpha_{j is})\, (s=1,2,\dots, n_{j i})$ in $Y^S$ 
are disjoint from the loop system $k_i\,(i=1,2,\dots,n)$,  
the 2-spheres $S(\alpha_{is})\,(s= 1,2,\dots, n_i)$ and 
$S(\alpha_{j is})\, (s=1,2,\dots, n_{j i})$ are considered in $S^4$. 
The 2-sphere $S(\alpha_{i1})$ is identified with $K_i$ in $S^4$ for all 
$i\, (i=1,2,\dots,n)$. 
The following claim is shown. 

\phantom{x}

\noindent{\bf (A.3)} The 2-spheres $S(\alpha_{is})\,(i=1,2,\dots,n;\, s= 2,3,\dots, n_i)$ and $S(\alpha_{jis})\, (i,j=1,2,\dots,n, j\ne i;\, s=1,2,\dots, n_{j i})$
form a trivial $S^2$-link in $S^4$. 

\phantom{x}

By (A.1) and (A.3), the $S^2$-link $L=\cup_{i=1}^n K_i$ is shown to be a ribbon $S^2$-link in $S^4$.

\phantom{x}

\noindent{\bf Proof of (A.3).} 
The loops $k^S_t\, (t=1,2,\dots,n)$ in $S^4$ bound disjoint disks 
$D^S_t\, (i=1,2,\dots,n)$ in $S^4$. 
Hence the loop $k^S_t$ in $S^4$ is isotopic to a band sum $k'_t$ of some parallel links 
$P_t(m_i)\, (i=1,2,\dots,n)$ of the meridian loops $m_i\, (i=1,2,\dots,n)$ of $K_i\, (i=1,2,\dots,n)$ 
in $S^4$. 
For a parallel $k^{S+}_t$ of $k^S_t$ in $S^4$, let $D^{S+}_t$ be a move of $D^S_t$ 
with $\partial D^{S+}_t=k^{S+}_t$ in $S^4$ so that 
the disk $D^{S+}_t$ is disjoint from the annuli $A_i\, (i=1,2,\dots,n)$ and 
the disks $D_{ji}\, (i,j=1,2,\dots,n; j\ne i)$. The  2-spheres 
$S(\alpha_{is})\,(i=1,2,\dots,n;\, s=1, 2,\dots, n_i)$ and 
$S(\alpha_{jis})\, (i,j=1,2,\dots,n, j\ne i;\, s=1,2,\dots, n_{j i})$ may  be disjoint from the disk $D^{S+}_t$ in $S^4$. 
By passing through a thickening $D^{S+}_t\times I $ 
of the disk $D^{S+}_t$ for every $t(\ne i)$ in $S^4$, the annulus 
$A_i$ and the disk $D_{j i}$ in $Y^S$ extend respectively in $S^4$ to an annulus 
$\bar A_i$ with $\partial \bar A_i=(-k^S_i)\cup  m_i$ and a disk
$\bar D_{j i}$ with $\partial \bar D_{j i}=m_j$. 
The annuli $\bar A_i\, (i=1,2,\dots,n)$ and the disks 
$\bar D_{j i}\, (i,j =1,2,\dots,n; j\ne i)$ should be disjoint in $S^4$. 
For $s\geq 2$, let $S(\partial \alpha_{is})$ be the two sphere union 
which is the boundary of a regular neighborhood $B(\partial\alpha_{is})$ of the two point set 
$\partial \alpha_{is}$ in $B(\alpha_{is})$. 
The 2-sphere $S(\alpha_{is})$ can be replaced by the 2-sphere obtained from 
$S(\partial \alpha_{is})$ by surgery along a 1-handle attaching to 
$S(\partial \alpha_{is})$ whose core is a subarc $\alpha'_{is}$ of $\alpha_{is}$ in 
$B(\alpha_{is})$. 
The following observation (whose proof is obvious) is used. 

\phantom{x} 

\noindent{\bf Observation~A.4} 
The 2-sphere $S'$ obtained from the 2-spheres 
$S^2\times\{0, 1\}$ by surgery along a 1-handle $h'$ thickening the arc 
$p\times[0,1]\, (p\in S^2)$ bounds the 
unique 3-ball $B'=\mbox{cl}(S^2\times [0,1]\setminus h')$. 
Further, let $S''$ obtained from the 2-spheres 
$S^2\times\{\frac{1}{4}, \frac{3}{4}\}$ by surgery along a 1-handle $h''$ thickening the arc $p\times[\frac{1}{4}, \frac{3}{4}]$, and 
$B''=\mbox{cl}(S^2\times[\frac{1}{4}, \frac{3}{4}]\setminus h'')$ the 3-ball bounded by $S''$. 
If the 1-handle $h'$ is thinner than the 1-handle $h''$, then the 3-ball $B''$ is in 
the interior of the 3-ball $B'$. 

\phantom{x}

Assume that the arc $\alpha_{is}$ cuts an innermost disk $\delta$ from  
the annulus $\bar A_i$. 
Then the arc $\alpha'_{is}$ is $\partial$-relatively isotopic to 
an arc $J$ in $m_i$ through the disk $\delta$, so that  
the arc $\alpha'_{is}$ joining  the two sphere union $S(\partial \alpha_{is})$ is 
$\partial$-relatively isotopic to  an arc $J$ joining the boundary 
$(\partial J)\times K_i$ of a spherical shell 
$J\times K_i$ of  the circle bundle $\partial D^2\times K_i$ with $J\subset\partial D^2$ 
for  a normal disk bundle $D^2\times L$ in $S^4$.  
Thus, the 2-sphere $S(\alpha_{is})$ is isotopic to the boundary 2-sphere 
$\partial\Delta(\alpha_{is})$ of a 3-ball $\Delta(\alpha_{is})$ in the spherical shell 
$J\times K_i$ (see \cite{HK}). 
Note that the 3-ball $\Delta(\alpha_{is})$ does not meet the $S^2$-link $L$ although 
the trace of this isotopy may meet $L$ since the disk $\delta$ may meet $L$. 
By continuing this process, it is seen from Observation~A.4 that the 2-spheres $S(\alpha_{is})\, (s=2,3,\dots, n_i)$ are isotopic to the disjoint 
boundary 2-spheres $\partial\Delta(\alpha_{is})\, (s=2,3,\dots, n_i)$
of an inclusive 3-ball family $\Delta(\alpha_{is})\, (s=2,3,\dots, n_i)$ in 
$D^2\times K_i$, where an {\it inclusive 3-ball family} is a family of finite number 
of 3-balls such that any two members $B_1$ and $B_2$ have the property 
\[B_1\subset \mbox{Int}(B_2),\quad B_2\subset \mbox{Int}(B_1), \quad \mbox{or}\quad B_1\cap B_2=\emptyset.\]
For the disk $\bar D_{j i}$, the same argument above can be applied to see that 
the 2-spheres $S(\alpha_{j is})\, (s=1,2,\dots, n_{j i})$ are 
isotopic to the disjoint boundary 2-spheres 
$\partial\Delta(\alpha_{j is})\, (s=1,2,\dots, n_{j i})$
of an inclusive 3-ball family $\Delta(\alpha_{j is})\, (s=1,2,\dots, n_{j i})$ in 
$D^2\times K_j$ with $j\ne i$. 
Thus, for every $i$, the 2-spheres $S(\alpha_{is})\,( s= 2,3,\dots, n_i)$ and 
$S(\alpha_{j is})\, (s=1,2,\dots, n_{j i})$ form a trivial $S^2$-link in $S^4$. 
Since the annuli $\bar A_i\, (i=1,2,\dots,n)$ and the disks 
$\bar D_{j i}\, (i,j =1,2,\dots,n; j\ne i)$ are disjoint, it can be seen that
the 2-spheres $S(\alpha_{is})\,(i=1,2,\dots,n;\, s= 2,3,\dots, n_i)$ and 
$S(\alpha_{j is})\, (i,j=1,2,\dots,n, j\ne i;\, s=1,2,\dots, n_{j i})$
form a trivial $S^2$-link in $S^4$ by varying 
the radius of the disk $D$ of the normal disk bundle $D\times L$ of $L$ 
for every $i$. 
This completes the proof of (A.3). $\square$

\phantom{x} 

This completes the proof of Free Ribbon Lemma. $\square$

\phantom{x}

\noindent{\bf Acknowledgments.} This work was partly supported by JSPS KAKENHI Grant Numbers JP19H01788, JP21H00978 and MEXT Promotion of Distinctive Joint Research Center Program JPMXP0723833165.

\phantom{x}

\end{document}